\documentclass[11pt,reqno]{amsart}

\usepackage{amsthm,amsfonts,amsmath,amscd,amssymb,epsfig}
\usepackage{enumerate,array}
\usepackage{hyperref}
\usepackage{graphicx}
\usepackage{color}
\usepackage{amsfonts,amssymb}
\usepackage{mathrsfs}
\usepackage{a4wide}

\theoremstyle{plain}

\newtheorem{neu}{}[section]
\newtheorem{Cor}[neu]{Corollary}
\newtheorem*{Cor*}{Corollary}
\newtheorem{Thm}[neu]{Theorem}
\newtheorem*{Thm*}{Theorem}
\newtheorem{Prop}[neu]{Proposition}
\newtheorem*{Prop*}{Proposition}
\theoremstyle{definition}
\newtheorem{Lemma}[neu]{Lemma}
\newtheorem*{Rmk*}{Remark}
\newtheorem{Rmk}[neu]{Remark}

\newtheorem*{Ex*}{Example}
\newtheorem{Qu}[neu]{Question}
\newtheorem*{Qu*}{Question}

\newtheorem{Def}[neu]{Definition}

\theoremstyle{remark}

\newtheorem*{note}{Note}

\theoremstyle{definition}

\newcommand{\p}{\partial}
\newcommand{\om}{\omega}

\newcommand{\pf}{\longrightarrow}

\newcommand{\N}{{\mathbb{N}}}

\newcommand{\R}{{\mathbb{R}}}

\newcommand{\M}{\mathcal{M}}
\newcommand{\LLL}{\mathscr{L}}

\newcommand{\Crit}{{\rm Crit}}

\newcommand{\Ham}{\mathrm{Ham}}

\newcommand{\LL}{\mathcal{L}}

\renewcommand{\AA}{\mathcal{A}}

\newcommand{\x}{\times}

\newcommand{\beq}{\begin{equation}}
\newcommand{\beqn}{\begin{equation}\nonumber}
\newcommand{\eeq}{\end{equation}}

\newcommand{\bea}{\begin{equation}\begin{aligned}}
\newcommand{\bean}{\begin{equation}\begin{aligned}\nonumber}
\newcommand{\eea}{\end{aligned}\end{equation}}

\numberwithin{equation}{section}

\begin{document}
\title[Existence of leafwise intersection points in the unrestricted case] {Existence of leafwise intersection points\\ in the unrestricted case}
\author{Jungsoo Kang}
\address{Department of Mathematics, Seoul National University,
Kwanakgu Shinrim, San56-1 Seoul, South Korea, Email:
hoho159@snu.ac.kr}
\begin{abstract}
In this article, we study the question of existence of leafwise intersection points for contact manifolds which are not necessarily of restricted contact type.
\end{abstract}

\thanks{ This work is supported by the Basic research fund 2009-0074356 funded by the Korean Government and Hi-Seoul Fellowship 2010.}
\maketitle

\section{Introduction}
The study of existence of leafwise intersection points has become an important aspect of Hamiltonian dynamics. Leafwise intersection points interpolate between periodic orbits and Lagrangian intersection points. To be more precise, when a coisotropic submanifold in a $2n$ dimensional symplectic manifold of has codimension $0$ resp. $n$, a leafwise intersection point coincides with a periodic orbit resp. a Lagrangian intersection point. In fact, the higher codimensional case has already been considered and explored in \cite{Gu,Ka}, yet in this paper we restrict our interests to the codimension one case.

Let $(\Sigma,\lambda)$ be a contact hypersurface in a symplectic manifold $(M,\om)$, that is $d\lambda=\om|_\Sigma$. The Hamiltonian vector field $X_F$ is defined implicitly by $i_{X_F}\om=dF$ for a time-dependent Hamiltonian function $F\in C^\infty(S^1\x M)$ where $dF$ is only the derivative with respect to $M$  and we call the time 1-map $\phi_F$ of its flow a Hamiltonian diffeomorphism. In addition, we denote by $\Ham(M,\om)$ the group of Hamiltonian diffeomorphisms defined on $(M,\om)$ and $\Ham_c(M,\om)$ the group of Hamiltonian diffeomorphisms generated by compactly supported Hamiltonian functions. The symplectic structure $\om$ determines the {\em characteristic line bundle} $\LL_\Sigma\subset T\Sigma$ over $\Sigma$:
$$\LL_\Sigma:=\{(x,\xi)\in T_x\Sigma\,|\,\om_x(\xi,\zeta)=0\textrm{ for all }\zeta\in T_x\Sigma\}.$$
$\Sigma$ is foliated by the leaves of the characteristic line bundle and we denote by $L_x$ the leaf through $x\in\Sigma$ of the characteristic foliation. We note that these leaves are spanned by the Reeb vector field $R$ of $\lambda$ which is characterized by $\lambda(R)=1$ and $i_R d\lambda_{|\Sigma}=0$. Then a {\em leafwise intersection point} of $\phi\in \Ham(M,\om)$ is by definition a point $x\in\Sigma$ such that $\phi(x)\in L_x$.

In this paper, we consider $\big(\Sigma\times(-1,\infty),d((r+1)\lambda)\big)$ the symplectization of $(\Sigma,\lambda)$ of dimension $2n-1$ where $r$ is the coordinate on $(-1,\infty)$.
\begin{Qu}
Given $\phi\in\Ham_c\big(\Sigma\x(\vartheta_1,\vartheta_2),d((r+1)\lambda)\big)$ for $-1<\vartheta_1<0<\vartheta_2<\infty$, does $\phi$ have a leafwise intersection point?
\end{Qu}
We give an affirmative answer to the question above for a class of Hofer small $\phi$ (see Section 2 for the definition of the Hofer norm $||\cdot||$) and symplectically fillable contact manifolds. Throughout this paper, we assume that $\Sigma\x(\vartheta_1,\vartheta_2)$ is symplectically embedded in a symplectic manifold $(M,\om)$ of dimension $2n$.

\begin{Def}
We denote by $\wp(\Sigma,\lambda)>0$ the {\em minimal period} of a Reeb orbit of $(\Sigma,\lambda)$ which is contractible in $M$. If there is no contractible Reeb orbit we set $\wp(\Sigma,\lambda)=\infty$.
\end{Def}

\begin{Def}
We call a symplectic manifold $(M,\om)$ {\em convex at infinity} if $(M,\om)$ is symplectomorphic to the positive part of the symplectization of a compact contact manifold at infinity. Furthermore, $(M,\om)$ is called {\em symplectically aspherical} if one has the equality $\om|_{\pi_2(M)}=0$.
\end{Def}

\noindent\textbf{Theorem A.} Let $(M,\om)$ be closed (or convex at infinity) and symplectically aspherical and $(\Sigma,\lambda)$ be a contact hypersurface with  $\Sigma\x(-\vartheta_1,\vartheta_2)$ being symplectically embedded in $M$. Moreover, we assume that $\phi=\phi_F$ for $F\in C^\infty\big(S^1\x\Sigma\x(\vartheta_1,\vartheta_2)\big)$ where $F$ is constant outside the region $\Sigma\x[\rho_1,\rho_2]$ for $\rho_1\in(\vartheta_1,0),\rho_2\in(0,\vartheta_2)$ and has Hofer norm $||F||\leq\wp(\Sigma,\lambda).$ Then $\phi_F\in\Ham_c(M,\om)$ has a leafwise intersection point.

\begin{Rmk}
A contact hypersurface $(\Sigma,\lambda)$ in a symplectic manifold $(M,\om)$ is of restricted contact type if the contact 1-form $\lambda$ is defined on the whole symplectic manifold $M$ and $d\lambda=\om$. The main difference of this paper with other results in \cite{AF1}, \cite{Gi} and \cite{Gu} is that we drop the condition of restricted contact type. Thus, our ambient symplectic manifold need not be exact and therefore can be closed. Moreover we remove the condition needed in \cite{AF1} that $\Sigma$ bounds a compact region in $M$ by means of the argument developed in \cite{Ka}.
\end{Rmk}

On the other hand, for a special perturbation $F$, we are able to find a leafwise intersection point in a symplectization of $\Sigma$ even though the symplectization is not convex at infinity.
\begin{Def}
We denote the support of Hamiltonian vector field $X_F$ by
\beq
\mathrm{Supp}X_F:=\mathrm{cl}{\{(x,r)\in\Sigma\x(-1,\infty)\,|\,X_F(t,x,r)\ne0 \textrm{ for some }t\in S^1\}}.
\eeq
\end{Def}

\noindent\textbf{Theorem B.}
Let $\phi\in \Ham\big(\Sigma\x(-1,\infty),d((r+1)\lambda)\big)$ be of the form $\phi=\phi_F$ for some $F\in C^\infty(S^1\x\Sigma\x(-1,\infty))$ where $X_F$ is generated by the Reeb vector field $R$ and the Liouville vecor field $Y$ (defined before Proposition \ref{prop:Liouville v.f.}).
If $\mathrm{Supp}X_F$ is compact and $||F||\leq \wp(\Sigma,\lambda)$, then $\phi$ has a leafwise intersection point.

\begin{Rmk}
If the Weinstein conjecture holds, we can show Theorem B in an easier way. It is reduced to find a self intersection point of $S^1$ in $S^1\x(-1,\infty)$ where $S^1$ is diffeomorphic to the Reeb orbit since there exist at least one periodic Reeb orbit. Therefore Theorem B follows with the assumption $||\phi||\leq e(S^1)$ where $e(S^1)$ is the {\em displacement energy} of $S^1$ in $\mathbb{R}^2$.
\end{Rmk}

\subsection{Idea of the proofs}
Leafwise intersection points arise as critical points of a perturbed Rabinowitz action functional. Therefore our proof is based on Rabinowitz Floer homology as in \cite{AF1}. The difference to \cite{AF1} is that we do not assume restricted contact type of $\Sigma$. In order to overcome this difficulty, we apply an auxiliary Rabinowitz action functional similar as in \cite{CFP}. However, in our situation, we have to perturb the auxiliary Rabinowitz action functional. Using the auxiliary Rabinowitz action functional we show that the moduli space of gradient flow lines of original Rabinowitz action functional can be compactified.
To prove Theorem B, we compare the difference of the two action functionals and examine the energy of holomorphic curves. And then we notice that for special $F$, gradient flow lines of Rabinowitz action functional remain in a tubular neighborhood, that means, gradient flow lines do not go to infinity and thus we do not need the condition of convex at infinity.

\subsection{History and related results}
The problem of existence of leafwise intersection points was addressed by Moser \cite{M}. Moser obtained existence results for simply connected $M$ and $C^1$-small $\phi$. Banyaga \cite{B} removed the assumption of simply connectedness. Hofer \cite{H} and Ekeland-Hofer \cite{EH} replaced the assumption of $C^1$-smallness by boundedness of the Hofer norm below a certain symplectic capacity for restricted contact type in $\mathbb{R}^{2n}$. Ginzburg \cite{Gi} extended the Ekeland-Hofer results to subcritical Stein manifolds. Dragnev \cite{D} obtained the result on the leafwise intersection problem to closed contact type submanifold in $\R^{2n}$. Albers-Frauenfelder \cite{AF1} proved the existence of leafwise intersection points for a restricted contact hypersurface whenever a Hamiltonian diffeomorphism satisfies the same Hofer smallness assumption as in this article.
By a different approach G\"urel \cite{Gu} also proved existence for higher codimensional restricted contact type under the Hofer smallness. Ziltener \cite{Z} also studied the question in a different way and obtained a lower bound of the number of leafwise fixed points under the assumption that the characteristic foliation is a fibration.

\subsubsection*{Acknowledgments}
The author is grateful to Urs Frauenfelder for many valuable discussions and suggestions. He also thanks Peter Albers and the anonymous referee for helpful remarks and comments on an earlier version of this paper.

\section{Preliminaries}
Let $(\Sigma,\lambda)$ be a contact hypersurface in a closed (or convex at infinity) and symplectically aspherical symplectic manifold $(M,\om)$. In addition, we consider a non-autonomous Hamiltonian function $F\in C^\infty\big(S^1\x\Sigma\x(\vartheta_1,\vartheta_2)\big)$ which is constant outside the region $\Sigma\x[\rho_1,\rho_2]$ for $-1<\vartheta_1<\rho_1<0<\rho_2<\vartheta_2<\infty$. We extend $F$ locally constant to the whole symplectic manifold $M$. The existence of leafwise intersection points of the extension of $F$ guarantees existence of leafwise intersection points of $F$ because the Hamiltonian vector field of the extension vanishes on $\Sigma\x\big((\vartheta_1,\rho_1)\cup(\rho_2,\vartheta_2)\big)$. For simplicity we denote the extension of $F$ by $F$ again. Moreover we additionally assume that $\Sigma$ bounds a compact region in $M$. In other words, $M-\Sigma$ consists of two component. This additional assumption will be removed in Step 4 in the proof of Theorem A by using the argument developed in \cite{Ka}.

We introduce the Liouville vector field $Y$ for $\Sigma$, that is $i_Y\om=\lambda$ on $\Sigma\x(\vartheta_1,\vartheta_2)$. Then $Y$ has the following well-known properties.

\begin{Prop}\label{prop:Liouville v.f.}
\cite{HZ} Let the vector field $Y$ be the Liouville vector field for $(\Sigma,\lambda)$. Then $Y$ satisfies
\beq
L_Y\om=\om \qquad\textrm{and}\qquad Y\pitchfork T\Sigma.
\eeq
\end{Prop}
\hfill $\square$

We introduce a cutoff function $\varphi$ to extend $\lambda$ globally. Let $\varphi : \mathbb{R}\rightarrow\mathbb{R}$ and $supp\varphi\subset(\vartheta_1,\vartheta_2)$ such that $\varphi(r)=r+1$ for  $r\in[\rho_1,\rho_2]$ and $\varphi'(r)\leq 1+\kappa$ for all $r\in\mathbb{R}$ and for some $\kappa>0$ satisfying
\beq
\frac{1+\rho_1}{\rho_1-\vartheta_1}<1+\kappa.
\eeq
Then we have a global one form
\beq\label{eq:beta}
\beta(y):=\left\{\begin{array}{ll}
\varphi(r)\lambda(x)\qquad & y=(x,r)\in\Sigma\x(\vartheta_1,\vartheta_2),  \\[2ex]
0\qquad & y\in M-\big(\Sigma\x(\vartheta_1,\vartheta_2)\big).\end{array}\right.
\eeq

Let $\phi_Y^t$ be the flow of the Liouville vector field. Fix $\delta_1 > 0$ so that $\phi_Y^t|_\Sigma$ is defined for
$|t|\leq\delta_1$. Then we can define a function $\widehat G$ near $\Sigma$ by $\widehat G\big(\phi_Y^t(x)\big)=t$ for $x\in\Sigma.$
Let $U_\delta:=\{x\in M|\,\, |\widehat G(x)|<\delta\}$ for $0<\delta<\delta_1$ and choose $\delta_0<\delta_1$ satisfying $U_{\delta_0}\subset\Sigma\x[\rho_1,\rho_2]$. Then we can extend $\widehat G$ to $G\in C^{\infty}(M)$ to be defined on the ambient manifold M.

\beq
G:M\rightarrow\mathbb{R}\,\,\ by \,\,\,
  G = \left\{ \begin{array}{ll}
 \widehat G & \textrm {on  \,\,\,\,\,\,\,\,\,\,\,\,\,\,   $U_{{\delta_0}/2}$}\\
 locally\,\, constant & \textrm  {outside \,\,  $U_{\delta_0}$. } \end{array}\right.\;
\eeq
\\
\begin{note}
$X_G|_\Sigma=R$  and  $G^{-1}(0)=\Sigma.$
\end{note}


Consider a smooth function $\chi\in C^{\infty}(S^1,\mathbb{R})$ with $\int_0^1\chi(t)dt=1$ and $supp(\chi)\subset(\frac{1}{2},1)$. Using
$\chi$ and G, we can define a time-dependent Hamiltonian $H: S^1\x M\rightarrow\mathbb{R}$ by $H(t,x)=\chi(t)G(x)$. We also consider an almost complex structure $J$ on $M$ which is SFT-like; that is, it splits on $\Sigma\x(\vartheta_1,\vartheta_2)$ with respect to $T(\Sigma\x(\vartheta_1,\vartheta_2))=\ker\lambda\oplus(\ker d\lambda\oplus \R\frac{\partial}{\partial r})$ as $J|_{(\ker d\lambda\oplus \R\frac{\partial}{\partial r})}$ is a complex structure which interchanges the Reeb vector field $R$ with $\frac{\partial}{\partial r}$; strictly speaking, $JR=\frac{\p}{\p r}$ and $J\frac{\p}{\p r}=-R$.

\begin{Prop}\label{prop:om-dbeta}
For every $u\in TM$ the following inequality holds
\beq
d\beta(u,Ju)\leq(1+\kappa)\om(u,Ju)\;
\eeq
\end{Prop}
\begin{proof}
For $u\in T(\Sigma\x(\vartheta_1,\vartheta_2))$, we can write $u=u_1+u_2$ with respect to the decomposition $T(\Sigma\x(\vartheta_1,\vartheta_2))=\ker\lambda\oplus\big(\ker d\lambda\oplus \frac{\partial}{\partial r}\big).$\\
On $\Sigma\x(\vartheta_1,\vartheta_2)$, recall that we have chosen $\varphi$ as $\varphi(r)\leq r+1$ and $\varphi'(r)\leq 1+\kappa$ so that
\bea
d\beta(u,Ju)&\leq\varphi(r)d\lambda(u_1,Ju_1)+\varphi'(r)dr\wedge\lambda(u_2,Ju_2)\\
&\leq(r+1)d\lambda(u_1,Ju_1)+(1+\kappa)dr\wedge\lambda(u_2,Ju_2)\\
&\leq(1+\kappa)\om(u,Ju).
\eea
The last inequality follows from $\om=d((r+1)\lambda)=(r+1)d\lambda+dr\wedge\lambda$. Outside of $\Sigma\x(\vartheta_1,\vartheta_2)$, $d\beta$ vanishes but $\om(\cdot,J\cdot)$ is positive definite.
Therefore the proposition is proved.
\end{proof}

Next, we recall the definition of Hofer norm.
\begin{Def}\label{def:Ham}
Let $F\in C^\infty_c(S^1\x M,\mathbb{R})$ be a compactly supported Hamiltonian function. We set
\beq
||F||_+:=\int_0^1\max_{x\in M} F(t,x) dt\qquad||F||_-:=-\int_0^1\min_{x\in M} F(t,x) dt=||-F||_+
\eeq
and
\beq
||F||=||F||_++||F||_-.\;
\eeq
For $\phi\in \Ham_c(M,\om)$ the Hofer norm is
\beq
||\phi||=\inf\{||F||\mid \phi=\phi_F\}.\;
\eeq

\end{Def}

\begin{Lemma}\label{lemma:norms equivalent}
For all $\phi\in\Ham_c(M,\om)$
\beq
||\phi||=|||\phi|||:=\inf\{||F||\mid \phi=\phi_F,\;F(t,\cdot)=0\;\;\forall t\in[\tfrac12,1]\}\;.
\eeq
\end{Lemma}
\begin{proof}
To prove $||\phi||\geq|||\phi|||$, pick a smooth monotone increasing map $r:[0,1]\to[0,1]$ with $r(0)=0$ and $r(\tfrac12)=1$. For $F$ with $\phi_F=\phi$ we set $F^r(t,x):=r'(t)F(r(t),x)$. Then a direct computation shows $\phi_{F^r}=\phi_F$, $||F^r||=||F||$, and $F^r(t,x)=0$ for all $t\in[\tfrac12,1]$. The reverse inequality is obvious.
\end{proof}
Thanks to the previous Lemma \ref{lemma:norms equivalent} we only need to consider a Hamiltonian function $F$ with time support on $(0,\frac{1}{2})$ to prove the main theorems.

\section{Rabinowitz Action functionals}
\subsection{Critical points}

We denote by $\LLL\subset C^\infty(S^1,M)$ the component of contractible loops in $M$.
For Hamiltonian functions $H$ and $F$ defined so far, the perturbed Rabinowitz action functional $\AA_{F}^H(v,\eta):\LLL\times\mathbb{R}\longrightarrow\mathbb{R}$ is defined as follows:
\beq
\AA_{F}^H(v,\eta)=-\int_{D^2}\bar v^*\om-\int_0^1F(t,v(t))dt-\eta\int_0^1H(t,v(t))dt\;
\eeq
where $\bar v:D^2\rightarrow M$ is a filling disk of $v$. The symplectic asphericity condition implies that the value of the above action functional is independent of the choice of filling discs.

We also define the auxiliary Rabinowitz action functional
\beq
\widehat\AA_{F}^H(v,\eta):=-\int_{D^2}\bar v^*d\beta-\int_0^1F(t,v(t))dt-\eta\int_0^1H(t,v(t))dt\;.
\eeq
where $\beta$ has been defined in \eqref{eq:beta}.
Furthermore, we will use the difference of two action functionals:
\beq
\AA:=\widehat\AA_{F}^H-\AA_{F}^H=\int_{D^2}\bar v^*(\om-d\beta).
\eeq
Critical points $(v,\eta)\in\Crit\AA_{F}^H$ satisfy
\beq\label{eqn:critical point equation}\left.
\begin{aligned}
&\partial_tv=X_{F}(t,v)+\eta X_H(t,v)\\[1ex]
&\int_0^1H(t,v)dt=0
\end{aligned}
\;\;\right\}
\eeq

Albers-Frauenfelder observed that a critical point of $\AA_{F}^H$ gives rise to a leafwise intersection point.

\begin{Prop}\label{prop:critical point answers question}
\cite{AF1} Let $(v,\eta)\in\Crit\AA_{F}^H$. Then $x=v(0)$ satisfies $\phi_F(x)\in L_x$. Thus, $x$ is a leafwise intersection point.
\end{Prop}
\begin{proof}
Since $F(t,\cdot)$ vanishes for $t\in(\frac{1}{2},1)$, we compute for $t\geq\frac{1}{2}$,
\bea
\frac{d}{dt}G(v(t))&=dG(v(t))[\partial_t v]\\
&=dG(v(t))\bigr[\underbrace{X_F(t,v)}_{=0}+\eta \chi(t)X_{G}(v)\bigr]\\
&=0
\eea
Since $\int_0^1H(t,v(t))dt=0$ and $G(v(t))$ is constant for $t\geq1/2$, $v(t)\in H^{-1}(0)=G^{-1}(0)=\Sigma$ for $t\in[\frac{1}{2},1]$. On the other hand, $H$ has the time support on $(\frac{1}{2},1)$, $v$ solves the equation $\partial_t v=X_F(t,v)$ on $[0,\frac{1}{2}]$. Therefore $v(\frac{1}{2})=\phi_F^{1/2}(v(0))=\phi_F^1(v(0))$ since $F=0$ for $t\geq\frac{1}{2}$. For $t\in(\frac{1}{2},1)$, $\partial_t v=\eta X_{H}(t,v)$ implies $x=v(0)=v(1)\in L_{v(\frac{1}{2})}$. Thus we conclude that $x\in L_{\phi_F(x)}$, this is equivalent to $\phi_F(x)\in L_x$.
\end{proof}

\subsection{Gradient flow lines}

From now on, we allow $s$-dependence on F as follows:  \\$F_s(t,x) =F_-(t,x)$ for $s\leq-1$, for some $F_-$ and $F_s(t,x)=F_+(t,x)$ for $s\geq1$, for some $F_+$. Moreover $F_s(t,\cdot)=0$ for $t\in(\frac{1}{2},1)$. We also choose a family $J(s)$ of compatible almost complex structures on $M$ so that they still are SFT-like and $J(s)=J_-$ for $s\leq-1$, for some $J_-$ and $J(s)=J_+$ for $s\geq1$, for some $J_+$.

\begin{note}
The time supports of $H$ and $F_s$ are disjoint.
\end{note}

On the tangent space $T_{(v,\eta)}(\LLL\times\mathbb{R})\cong T_v\LLL\x\mathbb{R}$, we define metric $m$ as follows:
\beq
m_{(v,\eta)}\big((\hat v_1,\hat\eta_1),(\hat v_2,\hat\eta_2)\big):=\int_0^1\om_v(\hat v_1,J\hat v_2)dt+\hat\eta_1\hat\eta_2\;.
\eeq
We also define another bilinear form $\widehat m$ with $\beta$.
\beq
\widehat m_{(v,\eta)}\big((\hat v_1,\hat\eta_1),(\hat v_2,\hat\eta_2)\big):=\int_0^1d\beta_v(\hat v_1,J\hat v_2)dt+\hat\eta_1\hat\eta_2\;.
\eeq

\begin{Def}
A map $w\in(v,\eta)\in C^\infty(\mathbb{R},\LLL\times\mathbb{R})$ which solves
\beq\label{eqn:gradient flow line}
\partial_s w(s)+\nabla_m \AA_{F_s}^H(w(s))=0.\;
\eeq
is called a gradient flow line of $\AA_{F_s}^H$.
\end{Def}

According to Floer's interpretation, gradient flow equation \eqref{eqn:gradient flow line} needs to be interpreted as $v:\mathbb{R}\x S^1\to M$ and $\eta:\mathbb{R}\to\mathbb{R}$ solving

\beq\label{eqn:gradient flow equation}\left.
\begin{aligned}
&\partial_sv+J_s(v)\bigr(\partial_tv-\eta X_{H}(t,v)-X_{F_s}(t,v)\bigr)=0\\[1ex]
&\partial_s\eta-\int_0^1H(t,v)dt=0 .
\end{aligned}
\;\;\right\}
\eeq

\begin{Def}
The energy of a map $w\in C^{\infty}(\R,\LLL\times\R)$ is defined as
\beq
E(w):=\int_{-\infty}^\infty||\partial_s w||_m^2ds\;.
\eeq
\end{Def}

\begin{Lemma}\label{lemma:energy estimate for gradient lines}
Let $w$ be a gradient flow line of $\nabla_m\AA_{F_s}^H$. Then
\beq\label{eqn:energy estimate for gradient lines}
E(w)\leq\AA_{F_-}^H(w_-)-\AA_{F_+}^H(w_+)+\int_{-\infty}^\infty||\partial_sF_s||_-ds\;.
\eeq
where $w_\pm=\lim_{s\to\pm\infty}w(s)$ and the negative part of Hofer norm $||\cdot||_-$ has been defined in Definition \ref{def:Ham}. Moreover, equality hold if $\partial_sF_s=0$.
\end{Lemma}

\begin{proof}
It follows from the gradient flow equation \eqref{eqn:gradient flow line}.
\bea
E(w)&=\int_{-\infty}^{\infty}m\bigr(-\nabla_m\AA^H_{F_s}(w(s)),\p_s w(s)\bigr)ds\\
&=-\int_{-\infty}^\infty d\AA_{F_s}^H(w(s))(\partial_sw(s))ds\\[1ex]
&=-\int_{-\infty}^\infty \frac{d}{ds}\Big(\AA_{F_s}^H(w(s))\Big)ds+\int_{-\infty}^\infty \big(\partial_s\AA_{F_s}^H\big)(w(s))ds\\[1ex]
&=\AA_{F_-}^H(w_-)-\AA_{F_+}^H(w_+)-\int_{-\infty}^\infty\int_0^1\partial_sF_s(t,v)dtds\\[1ex]
&\leq\AA_{F_-}^H(w_-)-\AA_{F_+}^H(w_+)+\int_{-\infty}^\infty||\partial_sF_s||_-ds\;.
\eea
\end{proof}

\begin{Prop}\label{prop:d}
If $(v,\eta) \in \LLL \times \R$ and  $(\hat{v},\hat{\eta}) \in
T_{(v,\eta)}(\LLL \times \R)=\Gamma(S^1,v^*TM)\x\mathbb{R}$ then the following assertion holds.
\beq
 d\widehat\AA^H_F(v,\eta)(\hat v,\hat\eta) = \widehat m\Bigl(\nabla_m\AA^H_F(v,\eta),(\hat v,\hat\eta)\Bigr).\;
\eeq
\end{Prop}
\begin{proof}
It holds
\bea
 dH=i_{X_H}\om = i_{X_H}d\beta \qquad\textrm{and}\qquad
 dF=i_{X_F}\om = i_{X_F}d\beta.
\eea
We know that
\bea
&i_{X_F}d\beta=i_{X_F}(\varphi(r)d\lambda + \varphi'(r)dr\wedge\lambda),\\
&i_{X_F}\om=i_{X_F}d\big((r+1)\lambda\big)=i_{X_F}\big((r+1)d\lambda+dr\wedge\lambda\big).
\eea
On the region $[\rho_1,\rho_2]\times\Sigma, ~ \varphi(r)=r+1$ implies $ \om=d\beta$. Outside the region $[\rho_1,\rho_2]\times\Sigma,$ $i_{X_F}\om=0=i_{X_F}d\beta$ by assumption. The other equality that $i_{X_H}\om=i_{X_H}d\beta$ is analogous to the above since we have chosen $\delta_0$ so that $U_{\delta_0}\subset\Sigma\x[\rho_1,\rho_2]$.\\
Next, we note the formula of $\nabla_m\AA^H_F$:
\begin{displaymath}\label{eq:nabla}
\nabla_m\AA^H_F=\left(\begin{array}{c}J(v)\big(\p_tv-\eta X_{H}(t,v)-X_F(t,v)\big)\\[1ex]
-\int_0^1H(t,v)dt\end{array}\right).
\end{displaymath}
Now it directly follows
\bea
d\widehat\AA^H_F(v,\eta)(\hat v,\hat\eta)
&=\int_0^1 d\beta\big(\p_t v,\hat{v})-\int_0^1dF(t,v)(\hat v)dt-\eta\int_0^1dH(t,v)(\hat v)dt-\hat{\eta}\int_0^1H(t,v)dt \\
&=\int_0^1 d\beta\big(\p_t v,\hat{v})-\om\big(\eta{X_H(t,v)}+{X_F(t,v)},\hat{v}\big) dt-\hat{\eta}\int_0^1H(t,v)dt \\
&=\int_0^1 d\beta\big(\p_t v-\eta{X_H(t,v)}-{X_F(t,v)},\hat{v}\big) dt - \hat{\eta}\int_0^1H(t,v)dt  \\
&=\widehat m\Bigl(\nabla_m\AA^H_F(v,\eta),(\hat v,\hat\eta)\Bigr).
\eea
\end{proof}

\begin{Prop}\label{prop:bound on AA}
Let a gradient flow line $w=(v,\eta)$ of $\AA_{F_s}^H$ converge asymptotically to $w_\pm:=\lim_{s\to\pm\infty}w(s)$. Then the following inequality holds.
\beq
\AA(w_-)-\kappa E(w)\leq\AA(w(s))\leq\AA(w_+)+\kappa E(w).
\eeq
\end{Prop}
\begin{proof}
Using Proposition \ref{prop:om-dbeta} and Proposition \ref{prop:d},
\bea\label{eq:dA}
\frac{d}{ds}\AA(w)
=& \frac{d}{ds}\widehat\AA^H_{F_s}(w) - \frac{d}{ds}\AA^H_{F_s}(w)\\[1ex]
=& d\widehat\AA^H_{F_s}(w)(\partial_s w) - d\AA^H_{F_s}(w)(\partial_s w)+\partial_s\widehat\AA_{F_s}^H(w)-\partial_s\AA_{F_s}^H(w)\\
=& \widehat m\Big(\nabla_m\AA^H_F(v,\eta),\p_s w\Bigr)-m\Big(\nabla_m\AA^H_F(v,\eta),\p_s w\Big)+\int_0^1\partial_sF_sdt-\int_0^1\partial_sF_sdt\\
=& \int_0^1(d\beta-\om)(-\p_s v,J\p_s v)dt -\Big(\int_0^1H(t,v)dt\Big)^2+\Big(\int_0^1H(t,v)dt\Big)^2\\
=& \int_0^1(\om-d\beta)(\p_s v,J\p_s v)dt\\
\geq&-\int_0^1\kappa\om(\p_s v,J\p_s v)dt.
\eea

\noindent Integrate both sides of \eqref{eq:dA} with respect to $s$ from $-\infty$ to $s_0\in\mathbb{R}$, then we get
\bea\label{eq:1}
\AA\big(w(s_0)\big)-\AA(w_-)
&= \int_{-\infty}^{s_0}\frac{d}{ds}\AA\big(w(s)\big)ds\\
&\geq-\kappa\int_{-\infty}^{s_0}\int_0^1\om\big(\p_s v,J\p_s v\big)dtds\\
&\geq-\kappa E(w).
\eea
\\
On the other hand, integrate from $s_0$ to $+\infty$ and obtain
\bea\label{eq:2}
\AA\big(w(s_0)\big)-\AA(w_+)
&=-\int_{s_0}^{\infty}\frac{d}{ds}\AA\big(w(s)\big)ds\\
&\leq \kappa\int_{s_0}^{\infty}\int_0^1\om\big(\p_s v,J\p_s v\big)dtds\\
&\leq\kappa E(w).
\eea
Combine above two inequalities \eqref{eq:1} and \eqref{eq:2}, then the proposition follows immediately.
\end{proof}

\begin{Prop}\label{prop:uniform bound of cutoff action functional}
 $\widehat \AA_{F_s}^H$ has uniform bounds along gradient flow lines of $\AA_{F_s}^H$ in terms of the asymptotic data, that is the action values of $\AA^H_{F_s}$ and $\widehat\AA^H_{F_s}$ at $w_\pm$.
 \end{Prop}
 \begin{proof}
 At first, let us show that $\AA_{F_s}^H$ is uniformly bounded along a gradient flow line $w(s)$.
\bea
 0&\leq-\int_{s_1}^{s_2}d\AA_{F_s}^H(w(s))(\partial_s w)ds\\
 &=\AA_{F_{s_1}}^{H}(w(s_1))-\AA_{F_{s_2}}^{H}(w(s_2)) - \int_{s_1}^{s_2}\int_0^1\partial_s F_s(t,v)dtds\\
&\leq\AA_{F_{s_1}}^{H}(w(s_1))-\AA_{F_{s_2}}^{H}(w(s_2))+\int_{s_1}^{s_2}||\partial_s F_s||_-ds.\\
\eea
 From above inequality we obtain
\bea
&\AA_{F_{s_2}}^{H}(w(s_2))\leq\AA_{F_-}^{H}(w_-)+\int_{-\infty}^{\infty}||\partial_s F_s||_-ds\\
&\AA_{F_{s_1}}^{H}(w(s_1))\geq\AA_{F_+}^{H}(w_+)-\int_{-\infty}^{\infty}||\partial_s F_s||_-ds
\eea
Therefore for any $s_0\in\mathbb{R}$, it holds
\beq\label{eq:uniform bound on AA^H_F}
\big|\AA_{F_s}^H(w(s_0))\big|\leq \max\{\AA_{F_-}^H(w_-),-\AA_{F_+}^H(w_+)\}+\int_{-\infty}^{\infty}||\partial_sF_s||_-ds.
\eeq
By the definition of $\AA$, we know
\beq
\big|\widehat\AA_{F_s}^H(w(s))\big|\leq\big|\AA_{F_s}^H(w(s))\big|+\big|\AA(w(s))\big|,
\eeq
but both terms on the righthand side are uniformly bounded in terms of the asymptotic data, recall Lemma \ref{lemma:energy estimate for gradient lines} and Proposition \ref{prop:bound on AA}. Hence the proposition is proved.
 \end{proof}

\begin{Thm}\label{thm:compactness of moduli}
Let $\M$ be a moduli space of gradient flow lines of $\AA_{F_s}^H$ with uniform action bounds of $\AA^H_{F_s}$ and $\widehat\AA^H_{F_s}$ like \eqref{eq:assumption}. Then this moduli space is compact modulo breaking. More specifically, for a sequence $\{w_n\}_{n\in\mathbb{N}}$ in $\M$ and for every reparametrization sequence $\sigma_n\in\R$ the sequence $w_n(\cdot+\sigma_n)$ has a subsequence which converges in $C^\infty_\mathrm{loc}(\R,\LLL\x\R)$.
\end{Thm}
Moreover if $w_n$ $C^\infty_\mathrm{loc}$-converges to $v$, we know $E(v)\leq\limsup_{n\in\N}E(w_n)$ by the following calculation.
$$
E(v)=\int_{-\infty}^\infty||\p_sw||^2ds=\lim_{T\to\infty}\int_{-T}^T||\p_sw||^2ds\leq\lim_{T\to\infty}\limsup_{n\in\N}E(w_n)=\limsup_{n\in\N}E(w_n).
$$
This observation will be used later in \eqref{eq:contradictorty estimation}.
\begin{proof}
If we establish the following facts, the proof of the theorem follows from standard arguments in Floer theory. For a sequence of elements $\{w_n=(v_n,\eta_n)\}_{n\in\mathbb{N}}$ in $\M$, we have
\begin{enumerate}
\item a uniform $L^\infty$-bound on $v_n$,
\item a uniform $L^\infty$-bound on $\eta_n$,
\item a uniform $L^\infty$-bound on the derivatives of $v_n$.
\end{enumerate}
(1) follows from the assumption ``convex at infinity". Once the $L^\infty$-bound on $\eta_n$ is established, the $L^\infty$-bound on the derivatives of $v_n$ follows from bubbling-off analysis together with the symplectic asphericity of $(M,\om)$. Hence Theorem \ref {thm:bound on eta} finishes the proof.
\end{proof}

\begin{Lemma} \label{lemma:bound on eta when gradient is small}
For $(v,\eta)\in\LLL\x\R$, there exist $\epsilon > 0$ and $C > 0$ such that
\beq
\big|\big|\nabla\mathcal{A}^H_{F_s}(v,\eta)\big|\big|_m\leq \epsilon
\quad \Longrightarrow \quad |\eta| \leq C
\bigr(\big|\widehat{\mathcal{A}}^H_{F_s}(v,\eta)\big|+1\bigr).
\eeq
\end{Lemma}
\begin{proof} The proof of lemma proceeds in three steps.\\
\textbf{Step\, 1:} Assume that $v(t)$ lies in $U_{\delta}=\{x\in M\,|\, |G(x)|<\delta\}$ for all $t\in (\frac{1}{2},1)$ where $\delta=\min\{1,\delta_0/2\}$. Then there exists a constant $C_1>0$ such that
\beq
|\eta|\leq C_1\bigr(\big|\widehat\AA_{F_s}^H(v,\eta)\big|+\big|\big|\nabla_m\AA_{F_s}^H(v,\eta)\big|\big|_m+1\bigr).
\eeq

We estimate
\bea
\big|\widehat\AA_{F_s}^H(v,\eta)\big|
=&\bigg|\int_0^1 v^*\beta+\eta\int_0^1 H(t,v(t))dt+\int_0^1 {F_s}(t,v(t))dt\bigg|\\
=&\bigg|\eta\int_0^1 \beta(v)\big(X_H(t,v)\big) dt+\int_0^1 \beta(v)\big(X_{F_s}(t,v)\big)dt\\
& +\int_0^1 \beta(v)\big(\partial_t v-\eta X_H(t,v)- X_{F_s}(t,v)\big)dt +\eta\int_0^1 H(t,v(t))dt+\int_0^1 {F_s}(t,v(t))dt\bigg|\\
\geq&\bigg|\eta\int_0^1\chi(t)\beta(v)\big(R(v)\big)dt\bigg|-\bigg|\int_0^1 \beta(v)(\partial_t v-\eta X_H(t,v)- X_{F_s}(t,v))dt\bigg|\\
& -\bigg|\eta\int_0^1 H(t,v(t))dt\bigg|-C_{\delta,{F_s}}\\
\geq&|\eta|-\delta|\eta|-C_\delta||\partial_t v-\eta X_H(t,v)- X_{F_s}(t,v)||_{L^1}-C_{\delta,F}\\[0.5ex]
\geq&(1-\delta)|\eta|-C_\delta||\partial_t v-\eta X_H(t,v)- X_{F_s}(t,v)||_{L^2}-C_{\delta,F}\\[0.5ex]
\geq&(1-\delta)|\eta|-C_\delta||\nabla_m\AA_{F_s}^H(v,\eta)||_{m}-C_{\delta,F}
\eea
where $C_\delta:= ||\beta|_{U_\delta}||_{L^\infty}$, $C_{\delta,F}:= ||F||_{L^\infty}+C_\delta||X_F||_{L^\infty}$ and $L^1$-, $L^2$-norms on $T\LLL$ are taken with respect to the metric $g(\cdot,\cdot)=\om(\cdot,J\cdot)$.\\
Thus we get
\beq
 |\eta|\leq\frac{1}{1-\delta}\bigg(\big|\widehat\AA_{F_s}^H(v,\eta)\big|+C_\delta\big|\big|\nabla_m\AA_{F_s}^H(v,\eta)\big|\big|_{m}+C_{\delta,F}\bigg).
\eeq
This proves Step 1 with
 \beq
 C_1:=\max\bigg\{\frac{1}{1-\delta},\frac{C_\delta}{1-\delta},\frac{C_{\delta,F}}{1-\delta}\bigg\}.
 \eeq
 \\
\textbf{Step\, 2:} There exists $\epsilon>0$ with $||\nabla_m \AA_{F_s}^H(v,\eta)||_m\geq\epsilon$ if there is $t\in (\frac{1}{2},1)$ such that $v(t)\notin U_\delta$.\\

 If $v(t) \in M-U_{\delta/2}$ for all $t \in (\frac{1}{2},1)$ then easily we have
\beq
||\nabla_m\AA_{F_s}^H(v,\eta)||_m\geq\bigg|\int_0^1H(t,v(t))dt\bigg|\geq\frac{\delta}{2}.\;
\eeq
Otherwise there exists $t'\in (\frac{1}{2},1)$ such that $v(t')\in U_{\delta/2}$. Thus we can find $t_0, t_1 \in (\frac{1}{2},1)$ such that
\bea
v(t_0)\in \partial U_{\delta/2},\,\, v(t_1)\in \partial U_{\delta},\qquad\mathrm{and}\qquad\forall s\in [t_0,t_1],\,\, v(s)\in U_\delta-U_{\delta/2}\
\eea
or
\bea
 v(t_1)\in \partial U_{\delta},\,\, v(t_0)\in \partial U_{\delta/2},\qquad\mathrm{and}\qquad\forall s\in [t_1,t_0],\,\, v(s)\in U_\delta-U_{\delta/2}.
\eea
 We treat only the first case; the later case is analogous. For $\mathfrak{G}:=\max_{x\in U_{\delta}}||\nabla G(x)||_g$, we have
\bea
\mathfrak{G}||\nabla_m\AA_{F_s}^H(v,\eta)||_m &\geq \mathfrak{G}||\partial_t v-\eta X_H(v)- X_{F_s}(v)||_{L^2}\\
&\geq \mathfrak{G}||\partial_t v-\eta X_H(t,v)- X_{F_s}(t,v)||_{L^1}\\
&\geq \int_{t_0}^{t_1}||\partial_t v-\eta X_H(t,v)- X_{F_s}(t,v)||_g||\nabla G(x)||_g dt\\
&\geq \bigg|\int_{t_0}^{t_1}\big\langle\nabla G(v(t)),\partial_t v-\eta X_H(t,v) - X_{F_s}(t,v)\big\rangle_g dt\bigg|\\
&= \bigg|\int_{t_0}^{t_1}dG(v(t))\bigr(\partial_t v-\eta X_H(t,v) - \underbrace{X_{F_s}(t,v)}_{=0}\bigr) dt\bigg|\\
&= \bigg|\int_{t_0}^{t_1}\frac{d}{dt} G(v(t))-\eta \underbrace{dG(v)\big(X_H(t,v)\big)}_{=\chi\om(X_G,X_G)=0}dt\bigg|\\
&\geq \big|G(v(t_1))\big|-\big|G(v(t_0))\big|\\
&=\frac{\delta}{2}.
\eea
Therefore Step 2 follows with $\epsilon=\min\bigr\{\frac{\delta}{2},\frac{\delta}{2\mathfrak{G}}\bigr\}.$\\
\\
\textbf{Step\, 3:} Proof of the lemma.\\
\\
Step 2 yields that $v(t)\in U_\delta$ for all $t\in(\frac{1}{2},1).$ Thus we are able to apply Step 1 and it shows that $|\eta| \leq C
\bigr(\big|\widehat{\mathcal{A}}^H_{F_s}(v,\eta)\big|+1\bigr)$ with $C=C_1+\epsilon.$
\end{proof}

For a gradient flow line $w$ of $\AA^H_{F_s}$ and $\sigma\in\mathbb{R}$, we set
\beq
\tau(w,\sigma):=\inf\bigr\{\tau\geq0\,\bigr|\,\,||\nabla_m\AA_{F_s}^H(w(\sigma+\tau))||_m\leq\epsilon\bigr\}.
\eeq

\begin{Lemma}\label{lemma:tau(sigma)}
We have a bound on $\tau(w,\sigma)$ as follows:
\beq
\tau(w,\sigma)\leq\frac{\AA_{F_-}^H(w_-)-\AA_{F_+}^H(w_+)+C_F}{\epsilon^2}
\eeq
for $C_F:=\int_{-\infty}^{\infty}||\p_s F_s||_-ds \,\,<\,\, \infty.$
\end{Lemma}
\begin{proof}
Using Lemma \ref{lemma:energy estimate for gradient lines}, we compute
\bea
\epsilon^2\tau(w,\sigma)
&\leq\int_\sigma^{\sigma+\tau(w,\sigma)}\big|\big|\nabla_m \AA_{F_s}^H(w)\big|\big|_m^2ds\\
&\leq E(w)\\
&\leq\AA^H_{F_-}(w_-)-\AA^H_{F_+}(w_+)+\int_{-\infty}^{\infty}||\p_s F_s||_-ds.
\eea
Dividing $\epsilon^2$ both sides, the lemma follows.
\end{proof}

\begin{Thm}\label{thm:bound on eta}
Assume that $w=(v,\eta)\in C^{\infty}(\mathbb{R},\LLL\x\mathbb{R})$ is a gradient flow line of $\AA_{F_s}^H$ for which there exist $a \leq b$ such that
\beq\label{eq:assumption}
\widehat\AA_{F_s}^H(w(s),\,\,\AA_{F_s}^H(w(s))\in[a,b] \quad\textrm{for all } s\in \mathbb{R}
\eeq
  Then the $L^\infty$-norm of $\eta$ is uniformly bounded.
\end{Thm}
\begin{proof}
Using the Lemma \ref{lemma:bound on eta when gradient is small} and Lemma \ref{lemma:tau(sigma)},
\bea
|\eta(\sigma)|&\leq|\eta(\sigma+\tau(w,\sigma))|+\int_\sigma^{\sigma+\tau(w,\sigma)}|\partial_s\eta(s)|ds\\
&\leq C\bigg(\Big|\widehat{\mathcal{A}}^H_{F_s}(v,\eta)\Big|+1\bigg)+\tau(w,\sigma)||H||_{L^\infty}\\
&\leq C\big(\max\{|a|,|b|\}+1\big)+\Big( \frac{|b-a|+C_F}{\epsilon^2}\Big)||H||_{L^\infty}.
\eea
\end{proof}
As we have already mentioned, Theorem \ref{thm:bound on eta} completes the proof of Theorem \ref{thm:compactness of moduli}.

\section{Proof of theorem A}
The proof of Theorem A proceeds in four steps. In first three steps we give a proof under the assumptions that $||F||<\wp(\Sigma,\lambda)$ and $\Sigma$ splits $M$ into two components. Step 4 finally removes these additional assumptions.\\[0.7ex]
\noindent\textbf{Step 1.} Theorem A holds when $||F||<\frac{\rho_1-\vartheta_1}{1+\rho_1}\wp(\Sigma,\lambda)$.
\begin{proof} For $0\leq r$, we choose a smooth family of functions $\beta_r\in C^{\infty}(\mathbb{R},[0,1])$ satisfying
\begin{enumerate}
\item for $r\geq1$: $\beta_r'(s)\cdot s\leq0$ for all $s\in\R$, $\beta_r(s)=1$ for $|s|\leq r-1$, and $\beta_r(s)=0$ for $|s|\geq r$,
\item for $r\leq1$: $\beta_r(s)\leq r$ for all $s\in\R$ and $\mathrm{supp}\beta_r\subset[-1,1]$,
\item $\lim_{r\to\infty}\beta_r(s\mp r)=:\beta_\infty^\pm(s)$ exists, where the limit is taken with respect to the $C^\infty_\mathrm{loc}$ topology.
\end{enumerate}

We fix a point $p\in\Sigma$ and consider the moduli space
\beq
\M:=\left\{(r,w)\in[0,\infty)\times C^{\infty}(\R,\LLL\times\R)\bigg|\; \begin{aligned} &w:\text{gradient flow line of $\AA^H_{\beta_rF},$}\\& w_-=(p,0),\;w_+\in\Sigma\x\{0\} \end{aligned}\right\}\;.
\eeq
\underline{Claim}: If there exists no leafwise intersection point, then $\M$ is compact. Moreover, its boundary consists of the point $(0,p,0)$ only.\\
Proof of Claim. For $(r,w)\in\M$,
\bea
E(w)&=-\int_{-\infty}^{\infty}d\AA_{\beta(s)_rF}^H(w(s))(\partial_s w) ds\\
&\leq\AA_0^H(p,0)-\AA_0^H(p,0)+\int_{-\infty}^\infty||\partial_s \beta_r(s)F||_-ds\\
&=\int_{-\infty}^\infty||\beta_r'(s)F||_-ds\\
&=\int_{-\infty}^0\beta_r'(s)||F||_-ds-\int_{0}^\infty\beta_r'(s)||F||_+ds\\
&=\beta_r(0)\big(||F||_-+||F||_+\big)\\
&\leq||F||.
\eea
Accordingly, we can also estimate
\beq\label{eq:step 1 in thm A}
-||F||<\AA^H_{\beta_r F}(w(s))<||F||,\qquad (r,w)\in\M.
\eeq
Since we have uniform action bounds, a sequence $\{(r_n,w_n)\}_{n\in\N}$ in $\M$ converges (after choosing a subsequence) up to breaking, see Theorem \ref{thm:compactness of moduli}. If breaking along a sequence occurs, then one of the following has to exist.\footnote{See the proof of Theorem A in \cite{AF1} for the detail arguments.}

\begin{enumerate}
\item a non-constant gradient flow line $v$ of $\AA^H_0$ with one asymptotic end being $(p,0)$,
\item a gradient flow line $v$ of $\AA_{\beta^{\pm}_{\infty}F}^H$, where $\beta^{\pm}_{\infty}$ is as above.
\end{enumerate}
However, we can rule out the second case since otherwise one asymptotic end of $v$ is a critical point of $\AA_F^H$ which gives a leafwise intersection point according to Proposition \ref{prop:critical point answers question} and it contradicts to our assumption. In the first case, not both asymptotic ends of $v$ can be $(q,0)$ for some $q\in\Sigma$, otherwise $E(v)=\big|\AA^H(p,0)-\AA^H(q,0)\big|=0$ and it yields that $v$ would be a constant curve. For this reason, one of the asymptotic ends of $v$ has to be of the form $(\gamma,\eta)$ where $\gamma$ is a Reeb orbit with nonzero period $\eta$. Let us assume that $(\gamma,\eta)$ is positive end of $v$, the other case is analogous. Then we have the following estimation; let $\{s_n\}_{n\in\mathbb{N}}$ be a sequence in $\R$ such that $w_n(s_n)$ converge to $v(+\infty)$, then we estimate

\bea\label{eq:contradictorty estimation}
|\eta|&=\big|\widehat\AA^H_0(\gamma,\eta)\big|\\
&\leq\big|\AA^H_0(\gamma,\eta)\big|+\big|\AA(\gamma,\eta)\big|\\
&=E(v)+\bigg|\int_{-\infty}^{\infty}\frac{d}{ds}\AA(v(s))ds\bigg|\\
&\leq\limsup_{n\in\mathbb{N}} E(w_n)+\limsup_{n\in\mathbb{N}}\bigg|\int_{-\infty}^{s_n} \frac{d}{ds}\AA(w_n(s))ds\bigg|\\
&\leq ||F||+\limsup_{n\in\mathbb{N}}\big|\sup_{s\in\mathbb{R}}\AA(w_n(s))-\AA(p,0)\big|\\
&\leq ||F||+\limsup_{n\in\mathbb{N}} \big(\max\{\big|\AA(w_{+})\big|,\big|\AA(w_{-})\big|\}+\kappa E(w)\big)\\
&\leq ||F||+\kappa||F||.
\eea
If $||F||<\frac{1}{1+\kappa}\wp(\Sigma,\lambda)$, then due to the previous estimation \eqref{eq:contradictorty estimation} we deduce the contradiction $|\eta|<\wp(\Sigma,\lambda)$. Since we have chosen any $\kappa$ satisfying $\frac{1+\rho_1}{\rho_1-\vartheta_1}<1+\kappa$, taking the limit $\kappa\to\frac{1+\rho_1}{\rho_1-\vartheta_1}-1$ we deduce a contradiction to the assumption $||F||<\frac{\rho_1-\vartheta_1}{1+\rho_1}\wp(\Sigma,\lambda)$. This proves the Claim.
\hfill $\square$\\
We are able to regard the moduli space $\M$ as the zero set of a Fredholm section of a Banach bundle over a Banach manifold. Moreover, the Fredholm section is already transversal at the boundary point since the boundary is a constant solution and at a constant solution $\AA^H$ is Morse-Bott. Since $\M$ is compact by the previous claim, we can perturb a Fredholm section away from the boundary point to get a transverse Fredholm section whose zero set is a compact manifold with a single boundary point $(0,p,0)$. But such a manifold does not exist. This finishes the proof of Step 1 by contradiction.
\end{proof}

\noindent\textbf{Step 2.} There exist a symplectic manifold $\widehat M$ which is symplectically aspherical and convex at infinity and a symplectic embedding $\psi':\Sigma\x(\vartheta'_1,\vartheta_2)\to \widehat M$ for any $-1<\vartheta'_1<\vartheta_1$.
\begin{proof}
We have additionally assumed that a closed contact manifold $\Sigma$ splits $M$ into two components. We call a bounded component $M_b$ with $\p M_b=\Sigma$. Then we get a new symplectic manifold $(\widehat M,\widehat\om)$ which is still symplectically aspherical and convex at infinity as follows:
\beq
\widehat M:=M_b\cup_{\p M_b}\Sigma\x[0,\infty)
\eeq
\beq
  \widehat\om = \left\{ \begin{array}{ll}
 \om & \textrm {on}\qquad M_b,\\[1ex]
 d((r+1)\lambda) & \textrm  {on}\qquad  \Sigma\x(\vartheta_1,\infty). \end{array}\right.\;
\eeq
$\widehat\om $ is well-defined since $\om$ equals to $d((r+1)\lambda)$ on $\Sigma\x(\vartheta_1,0]\subset M_1$. Then we rescale the symplectic structure $\widehat\om$. For arbitrary small $0<\nu<1$, we have a rescaled symplectic manifold $(\widehat M,\nu\widehat\om)$. Then there is a symplectic embedding
\bea
\psi_\nu:\bigr(\Sigma\x[-1+\nu,\infty),d((r+1)\lambda)\bigr)&\pf\bigr(\Sigma\x[0,\infty),d((r+1)\nu\lambda)\bigr)\subset
\bigr(\widehat M,\nu\widehat\om\bigr)\\
\big(x,r\big)&\longmapsto\big(x,\frac{1}{\nu}(r-\nu+1)\big).
\eea
Therefore $\Sigma\x(\vartheta'_1,\vartheta_2)$ can be embedded into $(\widehat M,\nu\widehat\om)$ via $\psi'=\psi_\nu$ for $0<\nu<1+\vartheta'_1$, and it finishes the proof of Step 2.
\end{proof}
\noindent\textbf{Step 3.} Proof of Theorem A for the case that $||F||<\wp(\Sigma,\lambda)$ and $\Sigma$ splits $M$ into two components.
\begin{proof}
Let the Hofer norm of $F$ be less than $\wp(\Sigma,\lambda)$. Thus we pick $\vartheta'_1>-1$ satisfying
\beq
||F||<\frac{\rho_1-\vartheta_1'}{1+\rho_1}\wp(\Sigma,\lambda).
\eeq
Then we can symplectically embed $\Sigma\x(\vartheta'_1,\vartheta_2)$ to the symplectic manifold $(\widehat M,\nu\widehat\om)$ with $0<\nu<1+\vartheta'_1$ by Step 2. Thus Step 1 enable us to find a leafwise intersection point.
\end{proof}
\noindent\textbf{Step 4.} End of the proof of Theorem A.
\begin{proof}
In the proof of Step 4, our contact hypersurface $\Sigma$ need not bound a compact region in $M$. We consider a family of time-dependent Hamiltonian functions $H_\nu\in C^\infty(S^1\x M)$ for $\nu\in\N$ where $H_\nu(t,x)=\chi(t)G_\nu(x)$ such that
\begin{enumerate}
\item $0<\epsilon_\nu<\min\{1,\delta_0/2\}$ converges to zero as $\nu$ goes to infinity,
\item for $x\in\Sigma$,
\beq
G_{\nu}|_{U_{2\epsilon_\nu}-U_{\epsilon_\nu/2}}(\phi^r_Y(x))= \left\{ \begin{array}{ll}
 r-\epsilon_\nu & \textrm{if}\,\, r>0\\[0.5ex]
 -r-\epsilon_\nu & \textrm{if}\,\, r<0, \end{array}\right.\;
\eeq
\item $G_{\nu}|_{M-U_{\delta_0}}=$ $constant$,
\item $G_\nu^{-1}(0)=\Sigma\x\{-\epsilon_\nu,\epsilon_\nu\}=:\Sigma_{-\epsilon_\nu}\cup\Sigma_{\epsilon_\nu}$.
\end{enumerate}
We note that $X_{G_{\nu}}|_{\Sigma_{\pm\epsilon_\nu}}=\pm R_{\pm\nu}$ where $R_{\pm\nu}$ is the Reeb vector field on $\Sigma_{\pm\epsilon_{\nu}}$, and we denote by $\phi^t_{R_{\pm\nu}}$ the flow of the Reeb vector field $R_{\pm\nu}$. Then according to Proposition \ref{prop:critical point answers question}, one of the followings holds: For $(v_\nu,\eta_\nu)\in\Crit\AA^{H_\nu}_F$,
\beq\label{eq:close leafwise intersections1}
\qquad\phi_F^1\big(v_\nu(\frac{1}{2})\big)=v_\nu(0)=\phi_{R_{+\nu}}^{-\eta_\nu}\big(v_\nu(\frac{1}{2})\big)\qquad\textrm{or}\\
\eeq
\beq\label{eq:close leafwise intersections2}
\phi_F^1\big(v_\nu(\frac{1}{2})\big)=v_\nu(0)=\phi_{R_{-\nu}}^{\eta_\nu}\big(v_\nu(\frac{1}{2})\big).
\eeq

Given a perturbation $F$ with $||F||<\wp(\Sigma,\lambda)$, the following holds for a sufficiently large $\nu\in\N$.
$$
||F||<\min\bigr\{\wp\big(\Sigma_{-\epsilon_\nu},(1-\epsilon_\nu)\lambda\big),\wp\big(\Sigma_{+\epsilon_\nu},(1+\epsilon_\nu)\lambda\big)\bigr\};
$$
then Step 1, 2, and 3 guarantee the existence of critical points $(v_\nu,\eta_\nu)$ of $\AA^{H_\nu}_F$. For clarity, let $\mathfrak{n}_\nu$ be $-\eta_\nu$ resp. $\eta_\nu$ and $\mathfrak{R}_\nu$ be $R_{+\nu}$ resp. $R_{-\nu}$ if \eqref{eq:close leafwise intersections1} resp. \eqref{eq:close leafwise intersections2} holds. Thus we have
\beq
\phi_F^1(v_\nu(\frac{1}{2}))=\phi_{\mathfrak{R}_{\nu}}^{\mathfrak{n}_\nu}(v_\nu(\frac{1}{2})).
\eeq
 Then estimation \eqref{eq:step 1 in thm A} in Step 1 implies the following lemma.

\begin{Lemma}\label{lemma:uniform bound on time}
$\mathfrak{n}_\nu$ is uniformly bounded in terms of $\lambda$ and $F$.
\end{Lemma}
\begin{proof}
We estimate like \eqref{eq:step 1 in thm A}.
\bea
||F||&\geq\big|\AA^{H_\nu}_F(v_\nu,\eta_\nu)\big|\\
&=\Big|\int_0^1v^*\lambda+\int_0^1 H_\nu(t,v_\nu(t))dt+\int_0^1F(t,v_\nu(t))dt\Big|\\
&=\Big|\int_0^1\chi(t)\lambda(v_\nu)\big(\pm\eta_{\nu} R_{\pm\nu}(v_\nu)+X_F(t,v_\nu)\big)dt+\int_0^1F(t,v_\nu(t))dt\Big|\\
&=\Big|\frac{\pm\eta_{\nu}}{1\pm\epsilon_\nu}+\int_0^1\lambda(v_\nu)\big(X_F(t,v_\nu)\big)+\int_0^1F(t,v_\nu(t))dt\Big|.\\
\eea
Therefore we conclude
\beq
|\mathfrak{n}_\nu|=|\eta_{\nu}|\leq2||F||+2||\lambda_{|U_{\delta_0/2}}||_{L^\infty}||X_F||_{L^\infty}+2||F||_{L^\infty}.
\eeq
\end{proof}
The two sequences of points $\{v_\nu(0)\}_{\nu\in\N}$ and $\{v_\nu(1/2)\}_{\nu\in\N}$ converge and we denote by
\bea
x_0:=\lim_{\nu\to\infty}v_\nu(0),\quad x_{1/2}:=\lim_{\nu\to\infty}v_\nu\big(\frac{1}{2}\big).
\eea
Obviously $x_0$ and $x_{1/2}$ are points in $\Sigma$. Moreover we know that
\beq\label{eq:limit 1}
x_0=\lim_{\nu\to\infty}v_\nu(0)=\lim_{\nu\to\infty}\phi_F^1(v_\nu(\frac{1}{2}))=\phi_F^1(\lim_{\nu\to\infty}v_\nu(\frac{1}{2}))=\phi_F^1(x_{1/2}).
\eeq
Furthermore, due to Lemma \ref{lemma:uniform bound on time}, we have a limit of $\{\mathfrak{n}_\nu\}_{\nu\in\N}$.
\beq
\lim_{\nu\to\infty}\mathfrak{n}_\nu=:\mathfrak{n}.
\eeq
Thus we conclude that $x_0$ and $x_{1/2}$ lie on a same leaf:
\beq\label{eq:limit 2}
x_0=\lim_{\nu\to\infty}v_\nu(0)=\lim_{\nu\to\infty}\phi_{\mathfrak{R}_{\nu}}^{\mathfrak{n}_\nu}(v_\nu(\frac{1}{2}))=\phi_{R}^{\mathfrak{n}}(x_{1/2}).
\eeq
It directly follows
\beq
\phi_{R}^{\mathfrak{n}}(x_{1/2})=\phi_F^1(x_{1/2})
\eeq
from equation \eqref{eq:limit 1} together with \eqref{eq:limit 2}.\\[0.5ex]

On the other hand, we consider a perturbation $F$ with $||F||=\wp(\Sigma,\lambda)$. Set $F_\mu:=\mu\cdot F$ for $\mu\in[0,1)$ then
$||F_\mu||<\wp(\Sigma,\lambda)$. By previous Step 1, 2, and 3 so far, we know the existence of a critical point of $\AA^{H_\nu}_{F_\mu}$,
namely $(v_{\nu,\mu},\eta_{\nu,\mu})\in\Crit\AA^{H_\nu}_{F_\mu}$. Using the same calculation in the proof of Lemma \ref{lemma:uniform bound on time},
we note that $\eta_{\nu,\mu}$ is uniformly bounded. Thus due to the Arzela-Ascoli's theorem, we are able to find $(v_\nu,\eta_\nu)\in\Crit\AA^{H_\nu}_F$
which gives rise to a leafwise intersection for $F$ in $\Sigma_{\pm\epsilon_\nu}$ and eventually we can find a leafwise intersection of $F$ in $\Sigma$
by applying the previous argument in Step 4. It completes the proof of Step 4, and hence the proof of Theorem A.
\end{proof}

\section{Proof of theorem B}

 From now on, we assume that $F_s$ is a time-dependent Hamiltonian function on the symplectization of $\Sigma$ such that $X_{F_s}$ is
spanned by the Reeb vector field $R$ and the Liouville vector field $Y$. We note that the Liouville vector field is of the form $(r+1)\frac{\p}{\p r}$ where $r$ is the coordinate on $(-1,\infty)$. Moreover, let $\mathrm{Supp}X_F$ be a compact subset of the symplectization as in Theorem B and $\pi:\mathrm{Supp}X_F\pf(-1,\infty)$ be a projection with respect to the $r$ coordinate, then we denote by
\beq
\varrho^-:=\min_{(x,r)\in\mathrm{Supp}X_F}\pi(x,r),\qquad \varrho^+:=\max_{(x,r)\in\mathrm{Supp}X_F}\pi(x,r).
\eeq
We are additionally able to define a Hamiltonian function $H$ appropriately on the symplectization as before since $\Sigma$ splits its symplectization and choose a different cut off function $\varphi : \mathbb{R}\rightarrow\mathbb{R}$ such that $supp\varphi\subset[\varrho^--\epsilon,\varrho^++\epsilon]$ for any small $\epsilon>0$ where $\varphi(r)=r-\varrho$ with $\varrho:=\min\{|\varrho^-|,|\varrho^+|\}$ on $[\varrho^-,\varrho^+]$ and $\varphi'(r)\leq 1$ for all $r\in\mathbb{R}$. Then we obtain a global one form $\beta=\varphi(r)\lambda$. Furthermore, we define action functionals $\AA_{F_s}^H$, $\widehat\AA_{F_s}^H$ and $\AA$ again with the new $\beta$ as before.

\begin{Prop}\label{prop:two assertions}
If $(v,\eta) \in \LLL \times \R$ and  $(\hat{v},\hat{\eta}) \in
T_{(v,\eta)}(\LLL \times \R)$ then the following two assertions
hold.
\begin{description}
 \item[(i)] $d\widehat\AA^H_{F_s}(v,\eta)(\hat v,\hat\eta) =
   \widehat m\Bigl(\nabla_m\AA^H_{F_s}(v,\eta),(\hat v,\hat\eta)\Bigr),$
 \item[(ii)] $
   (m-\widehat m)\Bigl((\hat v,\hat\eta),(\hat v,\hat\eta)\Bigr)
   \geq 0.$
\end{description}
\end{Prop}

\begin{proof}
Proofs of (i) and (ii) are almost the same as Proposition \ref{prop:d} and Proposition \ref{prop:om-dbeta} respectively. We have chosen $\varphi$ so that $\varphi'(r)=1$ on $\Sigma\x[\varrho^-,\varrho^+]$ and we have $\varphi(r)<r+1$, $\varphi'(r)\leq1$ for all $\mathbb{R}$ since $\varrho^->-1$. Thus the assertion (ii) follows from
\bea
d\beta(u,Ju)\leq\om(u,Ju)
\eea
using the computation in Proposition \ref{prop:om-dbeta}.
We have assumed that $X_F$ is spanned by $R$ and $Y$ so that
\beq
i_{X_F}d\beta=i_{X_F}(dr\wedge\lambda)=i_{X_F}\om.
\eeq
Therefore the assertion (i) follows from the exactly same argument in Proposition \ref{prop:d}.
\end{proof}

\begin{Cor}\label{cor:derivative of AA}
The action value of the functional $\AA=\widehat\AA_{F_s}^H-\AA_{F_s}^H$ is nondecreasing along a gradient flow line of $\AA_{F_s}^H$.
\end{Cor}
\begin{proof}
Using Proposition \ref{prop:two assertions}, we estimate with a gradient flow line $w(s)$ of $\AA_{F_s}^H$.
\bea
\frac{d}{ds}\AA(w(s))
&= \frac{d}{ds}\bigg(\widehat\AA_{F_s}^H(w(s))\bigg)-\frac{d}{ds}\bigg(\AA_{F_s}^H(w(s))\bigg)\\
&= d\widehat\AA_{F_s}^H(w)(\partial_s w)+(\partial_s \widehat\AA_{F_s}^H)(w)-d\AA_{F_s}^H(w)(\partial_s w)-(\partial_s \AA_{F_s}^H)(w)\\
&= m\bigr(\partial_s w,\partial_s w\bigr)-\widehat m\bigr(\partial_s w,\partial_s w\bigr)+\int_0^1\partial_s F_s(t,v) dt-\int_0^1\partial_s F_s(t,v) dt\\
&\geq\,0.
\eea
\end{proof}

\begin{Cor}\label{cor:AA=0}
$\AA(w(s))$ is identically zero for all $(r,w)\in\M$, the moduli space defined in the proof of Theorem A.
\end{Cor}
\begin{proof}
We note that $\AA(w_+)=\AA(w_-)=0$ since $w_\pm=w(\pm\infty)$ are constants. Therefore the proof immediately follows from the previous corollary.
\end{proof}


\begin{Prop}
Assume that $(r,w)=(r,v,\eta)$ is an element in $\M$. Then $v\in C^\infty(\R,\LLL)$ remains in $\Sigma\x[\varrho^-,\varrho^+]$.
 \end{Prop}
 \begin{proof}
Let us investigate the case that $v(s,t)$ goes out of the region $\Sigma\times[\varrho^--\epsilon,\varrho^++\epsilon]$.
Assume that $v(s,t)$ does not lie in $\Sigma\times[\varrho^--\epsilon,\varrho^++\epsilon]$ for $s_-<s<s_+$. It means that there exists a nonempty open subset $U\subset Z:=(s_-,s_+)\x S^1$ such that $v(s,t)\in\Sigma\x\big((-1,\varrho^--\epsilon)\cup(\varrho^++\epsilon,\infty)\big)$ for $(s,t)\in U$.\\
Using the previous corollary, we calculate
\bea
0=&\int_{s_-}^{s_+}\frac{d}{ds}\AA(w(s))\\
=&\int_{s_-}^{s_+}\int_0^1(\om-d\beta)(\partial_s v, J(v)\partial_s v) dtds\\
=&\int_{Z-U}(\om-d\beta)(\partial_s v, J(v)\partial_s v) dtds+\int_U\om(\partial_s v, J(v)\partial_s v) dtds.
\eea
The last equality holds since $d\beta$ vanishes on $\Sigma\x\big((-1,\varrho^--\epsilon)\cup(\varrho^++\epsilon,\infty)\big)$. However, $(\om-d\beta)(\partial_s v, J(v)\partial_s v)$ is bigger or equal to zero by the assertion (ii) in Proposition \ref{prop:two assertions} and $\int_U\om(\partial_s v, J(v)\partial_s v) dtds>0$. Thus this case cannot occur and accordingly every gradient flow line of $\AA_{\beta_r F}^H$ satisfying $w(-\infty)=(p,0)$ and $w(\-\infty)\in\Sigma$ lies in $\Sigma\times[\varrho^--\epsilon,\varrho^++\epsilon]$. Taking the limit $\epsilon\to 0$, this finishes the proof of the proposition.
\end{proof}

\begin{Rmk}
In the case that $v(s,t)$ for $s_-<s<s_+$ goes entirely out the region $\Sigma\times[\varrho^-,\varrho^+]$, we can show the above proposition more easily by using the energy argument. Since $d\beta$ and $H$ vanish outside of $\Sigma\x[\varrho^-,\varrho^+]$, we calculate
\bea
E(w)_{s_-}^{s_+}&=\int_{s_-}^{s_+}||\partial_s w||^2 ds\\
&=\int_{s_-}^{s_+}-\frac{d}{ds}\AA(w(s))\\
&=0.
\eea
It yields that $w$ is constant when $v$ is in the outside of the region, but such $w$ never exist.
\end{Rmk}
\quad\\
\noindent\textbf{Proof of Theorem B.} The previous proposition enable us to overcome the following problems, namely the $L^\infty$-bound on $v$ and the $L^\infty$-bound on the derivatives of $v$ although the symplectization of $\Sigma$ is not convex at infinity. The $L^\infty$-bound on $\eta$ is almost the same as what we showed and therefore Theorem \ref{thm:compactness of moduli} follows. Even easier, since $\om$ is exact on the symplectization, the bound follows from \cite{AF1}. Hence Theorem A guarantees the existence of a leafwise intersection point.
\hfill $\square$\\

\newpage


\end{document}